\documentclass{ifacconf}

\usepackage{natbib} 
\usepackage[centertags]{amsmath}
\usepackage{graphics}
\usepackage{amscd}
\usepackage{amsfonts}
\usepackage{amssymb}
\usepackage{helvet}
\usepackage{rotating}
\usepackage{epsfig}
\usepackage[latin1]{inputenc}
\usepackage{longtable}
\usepackage{natbib}
\usepackage[english,german,american]{babel}
\usepackage{relsize}
\usepackage[european]{circuitikz}
\usepackage{algorithmic}
\newcommand{\R}{\mathbb{R}}

\newcommand{\target}{T}
\newcommand{\servicer}{S}
\newcommand{\mass}{M}

\newcommand{\LVLHcoordinate}{\boldsymbol{\rho}}
\newcommand{\LVLHcoordinatex}{x}
\newcommand{\LVLHcoordinatey}{y}
\newcommand{\LVLHcoordinatez}{z}
\newcommand{\BFScoordinatex}{1}
\newcommand{\BFScoordinatey}{2}
\newcommand{\BFScoordinatez}{3}
\newcommand{\Qcoordinatei}{i}
\newcommand{\Qcoordinatej}{j}
\newcommand{\Qcoordinatek}{k}
\newcommand{\Qcoordinatel}{l}

\newcommand{\LVLHvelocityx}{\dot{x}}
\newcommand{\LVLHvelocityy}{\dot{y}}

\newcommand{\controlposition}{v}
\newcommand{\controlpositionvector}{\boldsymbol{\controlposition}}
\newcommand{\controlpositionx}{\controlposition_\LVLHcoordinatex}
\newcommand{\controlpositiony}{\controlposition_\LVLHcoordinatey}
\newcommand{\controlpositionz}{\controlposition_\LVLHcoordinatez}
\newcommand{\controlBFSpositionx}{\controlposition_\BFScoordinatex}
\newcommand{\controlBFSpositiony}{\controlposition_\BFScoordinatey}
\newcommand{\controlBFSpositionz}{\controlposition_\BFScoordinatez}
\newcommand{\quaternion}{q}
\newcommand{\quaternionvector}{\boldsymbol{\quaternion}}
\newcommand{\quaternionvectortarget}{\quaternionvector^{\target}}
\newcommand{\quaternionvectorservicer}{\quaternionvector^{\servicer}}
\newcommand{\quaternioni}{\quaternion_\Qcoordinatei}
\newcommand{\quaternionj}{\quaternion_\Qcoordinatej}
\newcommand{\quaternionk}{\quaternion_\Qcoordinatek}
\newcommand{\quaternionl}{\quaternion_\Qcoordinatel}
\newcommand{\quaterniondotx}{\dot{\quaternion}_\Qcoordinatei}
\newcommand{\quaterniondoty}{\dot{\quaternion}_\Qcoordinatej}
\newcommand{\quaterniondotz}{\dot{\quaternion}_\Qcoordinatek}
\newcommand{\quaterniondotw}{\dot{\quaternion}_\Qcoordinatel}
\newcommand{\orientation}{\omega}
\newcommand{\orientationvector}{\boldsymbol{\orientation}}
\newcommand{\orientationvectortarget}{\orientationvector^{\target}}
\newcommand{\orientationvectorservicer}{\orientationvector^{\servicer}}
\newcommand{\orientationx}{\orientation_\BFScoordinatex}
\newcommand{\orientationy}{\orientation_\BFScoordinatey}
\newcommand{\orientationz}{\orientation_\BFScoordinatez}
\newcommand{\orientationxtarget}{\orientation_\BFScoordinatex^\target}
\newcommand{\orientationytarget}{\orientation_\BFScoordinatey^\target}
\newcommand{\orientationztarget}{\orientation_\BFScoordinatez^\target}
\newcommand{\orientationxservicer}{\orientation_\BFScoordinatex^\servicer}
\newcommand{\orientationyservicer}{\orientation_\BFScoordinatey^\servicer}
\newcommand{\orientationzservicer}{\orientation_\BFScoordinatez^\servicer}
\newcommand{\orientationxtargetdot}{\dot{\orientation}_\BFScoordinatex^\target}
\newcommand{\orientationytargetdot}{\dot{\orientation}_\BFScoordinatey^\target}
\newcommand{\orientationztargetdot}{\dot{\orientation}_\BFScoordinatez^\target}
\newcommand{\orientationxservicerdot}{\dot{\orientation}_\BFScoordinatex^\servicer}
\newcommand{\orientationyservicerdot}{\dot{\orientation}_\BFScoordinatey^\servicer}
\newcommand{\orientationzservicerdot}{\dot{\orientation}_\BFScoordinatez^\servicer}
\newcommand{\inertia}{J}
\newcommand{\inertiavector}{\boldsymbol{\inertia}}
\newcommand{\inertiaxxtarget}{\inertia_{\BFScoordinatex\BFScoordinatex}^\target}
\newcommand{\inertiayytarget}{\inertia_{\BFScoordinatey\BFScoordinatey}^\target}
\newcommand{\inertiazztarget}{\inertia_{\BFScoordinatez\BFScoordinatez}^\target}
\newcommand{\inertiaxxservicer}{\inertia_{\BFScoordinatex\BFScoordinatex}^\servicer}
\newcommand{\inertiayyservicer}{\inertia_{\BFScoordinatey\BFScoordinatey}^\servicer}
\newcommand{\inertiazzservicer}{\inertia_{\BFScoordinatez\BFScoordinatez}^\servicer}
\newcommand{\controlmomentum}{m}
\newcommand{\controlmomentumvector}{\boldsymbol{\controlmomentum}}
\newcommand{\controlmomentumx}{\controlmomentum_\BFScoordinatex}
\newcommand{\controlmomentumy}{\controlmomentum_\BFScoordinatey}
\newcommand{\controlmomentumz}{\controlmomentum_\BFScoordinatez}
\newcommand{\state}{\boldsymbol{x}}
\newcommand{\transformationBFStoLVLH}{R}
\newcommand{\control}{\boldsymbol{u}}
\newcommand{\dockingpoint}{d}
\newcommand{\dockingvector}{\boldsymbol{\dockingpoint}}
\newcommand{\dockingvectortarget}{\dockingvector^\target}
\newcommand{\dockingvectorservicer}{\dockingvector^\servicer}
\newcommand{\dockingvectordottarget}{\dot{\dockingvector}^\target}
\newcommand{\dockingvectordotservicer}{\dot{\dockingvector}^\servicer}
\newcommand{\dockingpointxtarget}{\dockingpoint_1^\target}
\newcommand{\dockingpointytarget}{\dockingpoint_2^\target}
\newcommand{\dockingpointztarget}{\dockingpoint_3^\target}
\newcommand{\dockingpointxservicer}{\dockingpoint_1^\servicer}
\newcommand{\dockingpointyservicer}{\dockingpoint_2^\servicer}
\newcommand{\dockingpointzservicer}{\dockingpoint_3^\servicer}
\newcommand{\dockingpointLVLHxtarget}{\dockingpoint_\LVLHcoordinatex^\target}
\newcommand{\dockingpointLVLHytarget}{\dockingpoint_\LVLHcoordinatey^\target}
\newcommand{\dockingpointLVLHztarget}{\dockingpoint_\LVLHcoordinatez^\target}
\newcommand{\dockingpointLVLHxservicer}{\dockingpoint_\LVLHcoordinatex^\servicer}
\newcommand{\dockingpointLVLHyservicer}{\dockingpoint_\LVLHcoordinatey^\servicer}
\newcommand{\dockingpointLVLHzservicer}{\dockingpoint_\LVLHcoordinatez^\servicer}
\newcommand{\dockingpointdotLVLHxtarget}{\dot{\dockingpoint}_\LVLHcoordinatex^\target}
\newcommand{\dockingpointdotLVLHytarget}{\dot{\dockingpoint}_\LVLHcoordinatey^\target}
\newcommand{\dockingpointdotLVLHztarget}{\dot{\dockingpoint}_\LVLHcoordinatez^\target}
\newcommand{\dockingpointdotLVLHxservicer}{\dot{\dockingpoint}_\LVLHcoordinatex^\servicer}
\newcommand{\dockingpointdotLVLHyservicer}{\dot{\dockingpoint}_\LVLHcoordinatez^\servicer}
\newcommand{\dockingpointdotLVLHzservicer}{\dot{\dockingpoint}_\LVLHcoordinatex^\servicer}
\newcommand{\finaltime}{t_f}
\newcommand{\targetset}{\mathbb{T}}
\newcommand{\servicerset}{\mathbb{S}}
\newcommand{\feasibleset}{\mathbb{F}}
\newcommand{\costfunction}{J}
\newcommand{\weight}{\mu}
\newcommand{\weighttime}{\weight_{\finaltime}}
\newcommand{\weightthrust}{\weight_{\controlpositionvector}}
\newcommand{\weightmomentum}{\weight_{\controlmomentumvector}}
\definecolor{johannes}{rgb}{.9,.1,.1}
\definecolor{juergen}{rgb}{.1,.4,.9}

\begin{document}

\begin{frontmatter}

\title{Optimal Rendezvous Path Planning to an Uncontrolled Tumbling Target} 
\thanks[footnoteinfo]{This work was partially funded by the German Federal Ministry of Education and Research (BMBF), grant no. 05M10WNA.}
\author[UniBW]{J.~Michael} 
\author[UniBT]{K.~Chudej}
\author[UniBW]{M.~Gerdts}
\author[UniBW]{J.~Pannek}
\address[UniBW]{Faculty of Aerospace Engineering, University of the Federal Armed Forces, 85577 Munich, Germany (johannes.michael@unibw.de, matthias.gerdts@unibw.de, juergen.pannek@unibw.de).}
\address[UniBT]{Chair of Mathematics in Engineering, University of Bayreuth, 95440 Bayreuth, Germany (kurt.chudej@uni-bayreuth.de).}

\begin{abstract}
: As the number of uncontrollable objects in low earth orbit is rising, the thread of collisions and thus the breakdown of working satellites becomes worth analyzing. Consequently, projects on removing objects from the important orbits are taken into account by the international space associations. This paper is about the modelling and optimal path planning of a docking maneuver to an uncontrollable tumbling target. After deriving the system dynamics, we introduce boundary conditions to ensure a safe and realizable maneuver and a general Bolza type cost functional to incorporate different optimization goals. In order to solve the resulting problem, we transform the dynamics to a set of differential algebraic equations which allow us to employ a direct optimization method while perserving the energy of the system. The concluding simulation results show the reliability and effectiveness of this approach.
\end{abstract}


\end{frontmatter}

\section{Introduction}
\label{Section:introduction}
Since the first artificial earth satellite Sputnik was launched in 1957, the number of earth surrounding objects is increasing continuously. As a consequence of more satellites being brought into orbit without removing old and broken ones, the thread of possible collisions is rising. Yet, not only the satellites themself, but also operational debris like rocket bodies or fuel tanks raise the risk of collisions, see e.g. \cite{Taylor2006}. Another problem is that --- according to \cite{NASA:Remediation} --- the number of dangerous elements with a size of 10 cm will increase even if no further satellites are launched. This issue is due to the outcome of collisions as observed in the 2009 satellite crash which caused a cloud of around 600 new fragments, see \cite{Guardian}. The \cite{NASA:Remediation} states that the removal of five objects with high mass and collision risk per year can stop the long term growth of space debris in low earth orbits. \\
For this task methods have been developed on how removing such objects is to be accomplished. One is to catch the object with a net and drag it out of orbit, see e.g. \cite{Telegraph}. Another is to perform a docking maneuver with a service satellite on the object of interest and carry it in coupled state to a safe area. \cite{Fehse2003} distinguishs two different methods how the coupling can be archieved, namely docking and berthing. During docking an inflexible clamp on the servicer is joined with some suitable point in the target. In contrast to this direct coupling a berthing maneuver is performed by bringing the service satellite into a so called berthing box near the target, grabbing it with a mounted manipulator arm and using this connection to establish the desired coupled state. \\
In this paper we will focus on the modelling and optimal control of a docking maneuver. We present an extension of \cite{ChudejMichaelPannek2012} torwards a  more complicated motion of the target. In contrast to \cite{Boyarko2011} we concentrate on the calculation of the optimal control with a direct approach using an SQP method. Furthermore, a differential algebraic formulation is used to cope with the necessity to normalize quaternions during integration. The position controls are applied in body fixed coordinates and the satellite is not a sphere but rotational symmetric. The derived model uses a generalized docking condition and does not rely on a discretization of the state trajectory.  \\
The paper is organized as follows: In Section \ref{Section:setup} we derive the dynamics of the two spacecrafts including relative coordinates for the rendezvous. Furthermore, we introduce docking conditions as well as state and control constraints to guarantee a coupling of the two satellites at the end of the maneuver and ensure feasible trajectories. In Section \ref{Section:discretization}, we present the resulting optimal control problem and give details on the implementation. Last, we show simulation results in Section \ref{Section:simulation results} and conclude the paper with an outlook to future work.

\section{Model Setup}
\label{Section:setup}
In the following, we first derive the model for two spacecrafts in an orbit around the earth with a small relative distance in relation to orbit height. One is called the target ($\target$) and is supposed to be uncontrolled, while the other, the servicer ($\servicer$), is fully actuated in position and attitude. The dynamics of the two spacecrafts are assumed to be independent, that is for example the thrusters of the servicer do not influence the dynamics of the target and no collisions occur. Therefore, two non-interacting subsystems can be used to describe the motion of the considered model. Based on the model, we introduce terminal and boundary conditions in order to ensure a successful docking upon termination of the maneuver.

\subsection{Coordinate Systems}
For the derivation of the rendezvous model different coordinate systems are required. The first system we use is the local-vertical-local-horizontal system (LVLH). Its $\LVLHcoordinatex$ coordinate is defined as the extension of the connection of the earth center to the satellite pointing out of orbit, $\LVLHcoordinatey$ is pointing in the direction of movement and $\LVLHcoordinatez$ completes the orthogonal tripod. Additionally, we define a body fixed coordinate system for each spacecraft. This second system is supposed to have its origin at the center of mass of each satellite and its axes are aligned with the principal axis. Hence, the rotational dynamics can be modelled as a system of first order differential equations. Here, we suppose the unrotated state to be such that the body fixed and the LVLH axis coincide. To distinguish between coordinate systems, we utilize subscripts $\LVLHcoordinatex$, $\LVLHcoordinatey$ and $\LVLHcoordinatez$ to denote variables in LVLH systems and subscripts $\BFScoordinatex, \BFScoordinatey$ and $\BFScoordinatez$ for variables in body fixed coordinates. Details and illustrations of the coordinate systems can be found in \cite{Alfriend2009}. Throughout the work, we omit the explicit time dependency of variables to shorten the notation.

\subsection{Relative orbit dynamics}
To describe a docking maneuver, we need to derive relative equations of motion. To this end, the position of the service satellite is expressed in LVLH coordinates of the target. Starting with Kepler's two body problem for each spacecraft and adding a control vector $\controlpositionvector = [\controlpositionx \;\controlpositiony \; \controlpositionz ]^\top$ to the equation of the servicer, one obtains a second order nonlinear differential equation for the relative dynamics in the LVLH--system. Assuming the chief orbit to be circular and applying a first order taylor expansion on the right side we end up with the so called Hill--Clohessy--Wilshire--Equations
  \begin{align}
  	 \label{eq:cw_eq} 
    \ddot{\LVLHcoordinatex} &= 2n \LVLHvelocityy + 3n^{2} \LVLHcoordinatex  + \frac{\controlpositionx}{\mass}, \displaybreak[0] \nonumber \\
  	\ddot{\LVLHcoordinatey} &= - 2n \LVLHvelocityx + \frac{\controlpositiony}{\mass}, \displaybreak[0] \\
  	\ddot{\LVLHcoordinatez} &= - n^{2}\LVLHcoordinatez + \frac{\controlpositionz}{\mass}. \nonumber
  \end{align}
In these equations $\mass$ denotes the mass of the satellite and $n$ is the mean motion, a constant calculated out of the gravitiy constant and the orbit height. Note that higher terms are neglectable in this model since the relative distance is small compared to the operating altitude. A detailed derivation of these equations can be found, e.g., in \cite{Alfriend2009}. \\
Here, the position control vector $\controlpositionvector$ can be regarded as the applied thrust in orbit fixed coordinates. Later, we transform $\controlpositionvector$ into body coordinates to obtain the forces applied by the thrusters mounted on the satellite.

\subsection{Orientation dynamics}
In addition to the relative position of the spacecrafts it is necessary to model their orientation in space. To avoid the appearance of the gimbal lock phenomenon --- as it would occur in an Euler angle representation due to singularities --- we use unit quaternions as parametrization. Besides the avoidence of singularities quaternions posses the property that they, in contrast to a description with angles, are a continuous representation of orientations. As we will only provide the necessary features of quaternions, we refer to \cite{Tewari2007} and \cite{Wertz1978} for additional information.

A quaternion is a 4-tupel $\quaternionvector = \left[\quaternioni \; \quaternionj \; \quaternionk \; \quaternionl \right]^{\top}$ representing an orientation referring to an unrotated reference coordinate system. As reference system we use the LVLH--system. For all multiples of a quaternion representing the same rotation we restrict it to be of length one, meaning $\left\| \quaternionvector \right\| = \sqrt{\quaternioni^2 + \quaternionj^2 + \quaternionk^2 + \quaternionl^2} = 1$.

To transform a vector from the rotated coordinates to the unrotated system one has to multiply it from the left with the matrix
\begin{small}
\begin{align*}
	\transformationBFStoLVLH = \begin{bmatrix}
		\quaternioni^2 - \quaternionj^2 - \quaternionk^2 + \quaternionl^2 & 2 ( \quaternioni \quaternionj - \quaternionk \quaternionl ) & 2 ( \quaternioni \quaternionk + \quaternionj \quaternionl ) \\
		2 ( \quaternioni \quaternionj + \quaternionk \quaternionl ) & - \quaternioni^2 + \quaternionj^2 - \quaternionk^2 + \quaternionl^2 & 2 ( \quaternionj \quaternionk - \quaternioni \quaternionl ) \\
		2 ( \quaternioni \quaternionk - \quaternionj \quaternionl ) & 2 ( \quaternionj \quaternionk + \quaternioni \quaternionl ) & - \quaternioni^2 - \quaternionj^2 + \quaternionk^2 + \quaternionl^2
	\end{bmatrix}.
\end{align*}
\end{small}%
We like to point out that, in contrast to literature, we define this matrix as the rotation from a rotated to an unrotated state. For further observations on time dependent quaternions, we need to include the vector of angular velocities $\orientationvector = \left[\orientationx \; \orientationy \; \orientationz \right]^{\top}$, whose elements represent the angular velocities around the body fixed coordinate axes. This relationship is obtained by the first derivative of $\quaternionvector$
\begin{equation*}
 \label{eq:quat_deriv}
 \dot{\quaternionvector} = \frac{1}{2}\begin{bmatrix} \orientationvector \\ 0 \end{bmatrix} \otimes \quaternionvector,
\end{equation*}
where $\orientationvector$ is extended to a pure quaternion, which means that the fourth element is zero, and $\otimes$ is the quaternion multiplication. The derivation can be found, e.g., in \cite{StevensLewis2003}. Componentwise we obtain the matrix vector product
\begin{align}
	\label{eq:dynamic:quaternion}
	\begin{bmatrix}
      \quaterniondotx^\alpha \\
      \quaterniondoty^\alpha \\
      \quaterniondotz^\alpha \\
      \quaterniondotw^\alpha 
    \end{bmatrix} &= \frac{1}{2}
    \begin{bmatrix}
      0 & \orientationz^\alpha & -\orientationy^\alpha & \orientationx^\alpha \\
      - \orientationz^\alpha & 0 & \orientationx^\alpha & \orientationy^\alpha \\
      \orientationy^\alpha & - \orientationx^\alpha & 0 & \orientationz^\alpha \\
      - \orientationx^\alpha & - \orientationy^\alpha & - \orientationz^\alpha & 0
    \end{bmatrix}
    \begin{bmatrix}
      \quaternioni^\alpha \\
      \quaternionj^\alpha \\
      \quaternionk^\alpha \\
      \quaternionl^\alpha
    \end{bmatrix} \;
    \alpha \in \{\target, \servicer\}.
\end{align}
Here, the superscript $\alpha$ denotes that these dynamics have to be evaluated for the target (T) and for the servicer (S).

Last, we incorporate the change of the angular velocity of the two rotating bodies. Therefore we use Euler's gyroscopic equation
\begin{align*}
  \dot{\inertiavector} \cdot \orientationvector + \inertiavector \cdot \dot{\orientationvector} + \orientationvector \times \left( \inertiavector \cdot \orientationvector \right) = \controlmomentumvector,
\end{align*}
see, e.g. \cite{Chobotov1991}. With the above mentioned assumption that the fixed axis and the principle axis of the body coincide, the inertia tensor $\inertiavector$ has diagonal form $\inertiavector = diag(\inertiavector_{\BFScoordinatex\BFScoordinatex}, \; \inertiavector_{\BFScoordinatey\BFScoordinatey}, \; \inertiavector_{\BFScoordinatez\BFScoordinatez})$. Assuming that the mass distribution is constant over the time, meaning $\dot{\inertiavector} = 0$, the gyroscopic equation can be solved for $\dot{\orientationvector}$ and delivers the dynamics of the angular velocities
\begin{align}
	\label{eq:gyro}
	\orientationxservicerdot &= \frac{1}{\inertiaxxservicer} \left( \orientationyservicer \orientationzservicer \left( \inertiayyservicer - \inertiazzservicer \right) + \controlmomentumx \right), \displaybreak[0] \nonumber \\
	\orientationyservicerdot &= \frac{1}{\inertiayyservicer} \left( \orientationxservicer \orientationzservicer \left( \inertiazzservicer - \inertiaxxservicer \right) + \controlmomentumy \right), \displaybreak[0]  \nonumber \\
	\orientationzservicerdot &= \frac{1}{\inertiazzservicer} \left( \orientationxservicer \orientationyservicer \left( \inertiaxxservicer - \inertiayyservicer \right) + \controlmomentumz \right), \displaybreak[0]  \nonumber \\
	\orientationxtargetdot &= \frac{1}{\inertiaxxtarget} \left( \orientationytarget \orientationztarget \left( \inertiayytarget - \inertiazztarget \right), \right) \displaybreak[0]  \\
	\orientationytargetdot &= \frac{1}{\inertiayytarget} \left( \orientationxtarget \orientationztarget \left( \inertiazztarget - \inertiaxxtarget \right), \right) \displaybreak[0] \nonumber \\
	\orientationztargetdot &= \frac{1}{\inertiazztarget} \left( \orientationxtarget \orientationytarget \left( \inertiaxxtarget - \inertiayytarget \right) \right). \nonumber
\end{align}
Within \eqref{eq:gyro}, the vector $\controlmomentumvector = \left[\controlmomentumx \; \controlmomentumy \; \controlmomentumz \right]^{\top}$ denotes the applied torque of the service satellite around its body fixed coordinate axes, for example generated by control--moment--gyroscopes. As the target is supposed to be uncontrolled there is no additional momentum force on its right side. \\
Combining equations \eqref{eq:cw_eq} -- \eqref{eq:gyro},  we obtain a system of twenty first order differential equations for the dynamics of the two spacecrafts. The state vector consists of the components
\begin{align*}
	\state &:= [\LVLHcoordinatex, \LVLHcoordinatey, \LVLHcoordinatez, \dot{\LVLHcoordinatex}, \dot{\LVLHcoordinatey}, \dot{\LVLHcoordinatez}, \orientationxservicer, \orientationyservicer, \orientationzservicer, \quaternioni^\servicer, \quaternionj^\servicer, \quaternionk^\servicer, \quaternionl^\servicer, \\
	& \qquad \orientationxtarget, \orientationytarget, \orientationztarget, \quaternioni^\target, \quaternionj^\target, \quaternionk^\target, \quaternionl^\target]^\top,
\end{align*}
and the six controls applied to the servicer are summarized in the control vector
\begin{equation*}
	\control := [\controlpositionx, \controlpositiony, \controlpositionz, \controlmomentumx, \controlmomentumy, \controlmomentumz]^\top.
\end{equation*} 

\subsection{Docking conditions}
To model a successful docking maneuver between two spacecrafts we introduce docking constraints. To this end, docking points $\dockingvectorservicer = \left[\dockingpointxservicer \; \dockingpointyservicer \; \dockingpointzservicer \right]^{\top}$ and $\dockingvectortarget = \left[\dockingpointxtarget \; \dockingpointytarget \; \dockingpointztarget \right]^{\top}$ are defined in body fixed coordinates, representing for example some hook on the servicer and the exaust funnel on the target, see also \cite{Alfriend2009} and \cite{Boyarko2010} where similar models are used to express relative position and velocity of points lying on a solid body in space. Now we design the docking conditions to ensure that these points coincide in position and velocity at the end of the maneuver. As the points are defined in different coordinates we transform both into the LVLH system and calculate the distance between these points. The transformation into the LVLH system is done by a multiplication from the left with the previously defined rotation matrix $\transformationBFStoLVLH$. The translation from the targets LVLH system to the servicers is the relative distance $\LVLHcoordinate = \left[\LVLHcoordinatex \; \LVLHcoordinatey \; \LVLHcoordinatez \right]^{\top}$. Hence, we obtain  $\left[\dockingpointLVLHxtarget \; \dockingpointLVLHytarget \; \dockingpointLVLHztarget \right]^{\top} := \transformationBFStoLVLH^{T}\dockingvectortarget$ and $\left[\dockingpointLVLHxservicer \; \dockingpointLVLHyservicer \; \dockingpointLVLHzservicer \right]^{\top}  := \LVLHcoordinate + \transformationBFStoLVLH^{S}\dockingvectorservicer$ which gives us
\begin{align*}
	\begin{bmatrix}\dockingpointLVLHxservicer \\ \dockingpointLVLHyservicer \\ \dockingpointLVLHzservicer \end{bmatrix} - \begin{bmatrix}\dockingpointLVLHxtarget \\ \dockingpointLVLHytarget \\ \dockingpointLVLHztarget \end{bmatrix} = \LVLHcoordinate + \transformationBFStoLVLH^{S}\dockingvectorservicer - \transformationBFStoLVLH^{T}\dockingvectortarget.
\end{align*}
The relative velocity in LVLH coordinates can be derived by taking the first derivative of the previous equation. As the bodies are supposed to be fixed it holds that $\dockingvectordotservicer = \dockingvectordottarget = 0$ and so the derivative simplifies to
\begin{align*}
	\begin{bmatrix}\dockingpointdotLVLHxservicer \\ \dockingpointdotLVLHyservicer \\ \dockingpointdotLVLHzservicer \end{bmatrix} - \begin{bmatrix}\dockingpointdotLVLHxtarget \\ \dockingpointdotLVLHytarget \\ \dockingpointdotLVLHztarget \end{bmatrix} = \dot{\LVLHcoordinate} + \transformationBFStoLVLH^{\servicer} \orientationvectorservicer \times  \transformationBFStoLVLH^{S}\dockingvectorservicer - \transformationBFStoLVLH^{\target} \orientationvectortarget \times \transformationBFStoLVLH^{T}\dockingvectortarget.
\end{align*}
Now, we define the docking conditions to be fulfilled upon termination of the maneuver by setting the relative distance and the relative velocity of the docking points to zero, i.e. 
\begin{align*}
\begin{bmatrix}\dockingpointLVLHxservicer(\finaltime) \\ \dockingpointLVLHyservicer(\finaltime) \\ \dockingpointLVLHzservicer(\finaltime) \end{bmatrix} - \begin{bmatrix}\dockingpointLVLHxtarget(\finaltime) \\ \dockingpointLVLHytarget(\finaltime) \\ \dockingpointLVLHztarget(\finaltime) \end{bmatrix} = 0\; \text{and} \;
\begin{bmatrix}\dockingpointdotLVLHxservicer(\finaltime) \\ \dockingpointdotLVLHyservicer(\finaltime) \\ \dockingpointdotLVLHzservicer(\finaltime) \end{bmatrix} - \begin{bmatrix}\dockingpointdotLVLHxtarget(\finaltime) \\ \dockingpointdotLVLHytarget(\finaltime) \\ \dockingpointdotLVLHztarget(\finaltime) \end{bmatrix} = 0.
\end{align*}

\subsection{State and control constraints}
\label{Subsection:constraints}
Complementing the docking constraints, we have to take into account state and control constraints resembling realistic behaviour of the model. The state constraints are necessary to prevent collisions of the objects. To this end, we define the two sets $\targetset(\state(t))$ and $\servicerset(\state(t))$ in $\R^3$ as the areas of occupied space by the target and the servicer respectively which depend on the current position and orientation of the objects. The shape of each set can be described as a polyhedron with the center of mass as reference point and which is then rotated and translated according to the current state. \\
Assuming that there are no other objects close to the servicer--target pair, we can define the set of feasible states via
\begin{align*}
   \feasibleset := \left\{\state(t) \;|\; \targetset(\state(t)) \,\cap\, \servicerset(\state(t)) = \emptyset \right\}.
\end{align*}
Here, we like to note that the docking point definition and the feasible set must not contradict each other. If this was the case, then the set of feasible terminal states would be empty rendering the docking problem unsolvable.\\
Considering the control inputs, we impose bounds on both the thrusts and the applied momentum. These restrictions are motivated by technical limitations of the devices, such as the maximal momentum of the control moment gyroscopes and thust produced by the mounted nozzles. As the Hill--Clohessy--Wilshire--Equations are formulated in LVLH coordinates the position control is also aligned to these axes. To obtain the actual thrust applied by the thrusters mounted on the servicer, we rotate the control according to its current orientation via
\begin{align*}
 \begin{bmatrix} \controlBFSpositionx \\ \controlBFSpositiony \\ \controlBFSpositionz \end{bmatrix} = 
 \transformationBFStoLVLH^{\top} \begin{bmatrix} \controlpositionx \\ \controlpositiony \\ \controlpositionz \end{bmatrix}.
\end{align*}
Utilizing the controls in body fixed coordinates, we define the respective bounds via
\begin{align*}
 \left\| \transformationBFStoLVLH^{\top} \controlpositionvector \right\|_{\infty} &\leq \controlposition_{max} \\
 \left\| \controlmomentumvector \right\|_{\infty} &\leq \controlmomentum_{max} .
\end{align*}
Note that the constraint on $\controlpositionvector$ is defined in the infinity norm $\left\| \cdot \right\|_{\infty}$ which is not invariant under rotation. Hence, the transformation to body--fixed coordinates cannot be omitted.

\subsection{Cost functional}
As optimization goal, we minimize the Bolza type cost functional
\begin{align*}
 \costfunction(\control, \finaltime) = \weighttime \finaltime + \weightthrust \left\| \controlpositionvector(t) \right\|_{L_2}^2 + \weightmomentum \left\| \controlmomentumvector(t) \right\|_{L_2}^2.
\end{align*}
Here, the final time $\finaltime$ and costs for thrust and momentum control are combined as a weighted sum over the maneuvering period $[0, \finaltime]$. With the nonnegative weights $\weighttime, \weightthrust$ and $\weightmomentum$ the optimization goal can be modified towards time--optimality or a minimum energy consumption strategy. A different weighting of position and momentum control is motivated by the fact that the control momentum gyroscopes run with electric energy and can be recharged during or after the maneuver while the thrusters rely on some propellant which is limited in amount.

\section{Discretization and implementation}
\label{Section:discretization}

Combining all elements from the previous section, we obtain the following optimal control problem:
\begin{align}
	\label{opt_prob}
	&\mbox{Minimize } \quad J(\control,\finaltime) \nonumber\\ 
	 &\text{subject to the dynamics} \quad \eqref{eq:cw_eq} - \eqref{eq:gyro} \nonumber\\
	 &\text{with initial and terminal conditions} \nonumber\\
	  & \qquad \boldsymbol{x}(0) = \boldsymbol{x}_0 \displaybreak[0] \nonumber\\
	  & \qquad \boldsymbol{0} = \LVLHcoordinate(\finaltime) + \transformationBFStoLVLH^{\servicer}(\finaltime) \dockingvectorservicer - \transformationBFStoLVLH^{\target}(\finaltime)\dockingvectortarget \displaybreak[0] \nonumber\\
	  & \qquad \boldsymbol{0} = \dot{\LVLHcoordinate}(\finaltime) + \transformationBFStoLVLH^{\servicer}(\finaltime) \orientationvector^{\servicer}(\finaltime) \times \transformationBFStoLVLH^{\servicer}(\finaltime) \dockingvectorservicer \\
	  & \qquad \qquad \qquad \, - \transformationBFStoLVLH^{\target}(\finaltime) \orientationvector^{\target}(\finaltime) \times \transformationBFStoLVLH^{\target}(\finaltime)\dockingvectortarget \displaybreak[0] \nonumber\\
	    & \text{and constraints} \nonumber\\
	  & \qquad \left\| \transformationBFStoLVLH(t)^{\top} \controlpositionvector(t) \right\|_{\infty} \leq \controlposition_{max} \quad \forall t \in [0, \finaltime] \displaybreak[0] \nonumber \\
	  & \qquad \left\| \controlmomentumvector(t) \right\|_{\infty} \leq m_{\max} \quad \forall t \in [0, \finaltime] \displaybreak[0] \nonumber \\
	  & \qquad \state(t) \in \feasibleset \quad \forall t \in [0, \finaltime]. \nonumber
\end{align}
To solve the above problem, we impose a recursive direct approach. In contrast to a full discretization in which both states and controls are treated as optimization variables and the dynamics represent additional constraints, only the discretized controls are considered to be optimization variables and the states remain untouched and are not included in the constraint set. Instead, the state vectors at the discretization time instances are computed outside of the optimization problem. Hence, the resulting optimization problem is smaller and denser compared to a full discretization. In particular, here we use an equidistant discretization of the time interval $[0, \finaltime]$ into $N$ subintervals and apply the control in a sample and zero--order--hold scheme, that is the control is kept constant within each subinterval and may jump at the bounds of the subintervals. Consequently, the number of optimization variables of $20 ( N + 1 ) + 6 N$ is reduced to $6 N$ for the recursive discretization. \\
The dynamics are solved via an implicit linearized Runge--Kutta method of order two, see also \cite{Gerdts2005} for further details. Yet, since the unitarity property of the quaternions may not be perserved by the integration method, we reformulate the dynamics as a set of differential algebraic equations. Note that this property is vital since the terminal conditions heavily depend on the rotation of the body fixed coordinate systems. Here, we transform the fourth element of the servicer and target quaternions into the algebraic variables $\lambda^{\target}$ and $\lambda^{\servicer}$ and add the algebraic equation
\begin{equation*}
 0 = 1 - \left((\quaternioni^{\alpha})^2 + (\quaternionj^{\alpha})^2 + (\quaternionk^{\alpha})^2 + (\lambda^{\alpha})\right) \quad \alpha \in \left\{\target, \servicer \right\}.
\end{equation*}
to the dynamics. Due to this modification, energy preservation and hence the unitarity property of the quaternions is guaranteed. In this context, we like to stress that in contrast to a full discretization approach, a recursive method does not cause an additional optimization related error within the solution of the dynamics. Hence, although a full discretization may be faster to solve due to the comparably simple parallelizability of the underlying multiple shooting method, respective results may fail to satisfy the physical properties of the dynamics. \\
Apart from the dynamics, we also convert the cost functional $J(\control,\finaltime)$ from a Bolza to a Mayer type functional by introducing two additional state variables resembling the consumed thrust and momentum control. The respective values are evaluated at $\finaltime$ and incorporated into the cost functional using the respective weights $\weightthrust$ and $\weightmomentum$. The free final time property of Problem \eqref{opt_prob} is treated by transforming it to a fixed time interval and adding the terminal time $\finaltime$ to the list of optimization parameters. \\
Assuming a rotational symmetric satellite with respect to its $\BFScoordinatey$ axis and analyzing equation \eqref{eq:gyro}, we additionally obtain that $\orientationy^{\target}$, remains unchanged during the process while the other elements interchange their energy. Hence, the docking point of the target is moving on a circle and thus depends on the final time $\finaltime$. 
The resulting optimization problem is solved using the software package OCPID-DAE1, cf. \cite{Gerdts2010}, which applies a robust SQP method combined with a gradient calculation using sensitivity DAEs and is suitable for optimal control problems subject to differential algebraic equations of index one.

\section{Simulation Results}
\label{Section:simulation results}
In order to illustrate our approach, we consider the case of a tumbling target, i.e. an object which has a deviation from a stable rotation around its $\LVLHcoordinatey$--axis. Due to this deviation the body fixed $\BFScoordinatey$--axis describes a tumbling motion around the LVLH $\LVLHcoordinatey$--axis. Such a motion can occur if a satellite runs out of energy to maintain its stable rotation or when it is hit by an object. At the beginning of the maneuver, the servicer is supposed to be in an unrotated position with a relative distance of $10 m$ behind its target revealing the following exemplary initial values:
\begin{align*}
\LVLHcoordinate_0 & = 
\begin{bmatrix}
0 &
-10 &
0 
\end{bmatrix}^\top &
\dot{\LVLHcoordinate}_0 & = 
\begin{bmatrix}
0 &
0 &
0 
\end{bmatrix}^\top \displaybreak[0] \\
\quaternionvectorservicer_0 & =
\begin{bmatrix}
0 &
0 &
0 &
1 
\end{bmatrix}^\top &
\quaternionvectortarget_0 & =
\begin{bmatrix}
-0.05 &
0 &
0 &
0.99875 
\end{bmatrix}^\top \displaybreak[0]  \\
\orientationvectortarget_0 & =
\begin{bmatrix}
0 &
0.0349 &
0.017453
\end{bmatrix}^\top &
\orientationvectorservicer_0 & =
\begin{bmatrix}
0 &
0 &
0 
\end{bmatrix}^\top
\end{align*}
We like to note that a docking maneuver cannot be performed for any arbitrary tumbling target satellite. Imagine a motion where the docking point describes a very quick change in position and/or velocity as it occurs for example when the docking vector is orthogonal to the main rotation axis. Yet, such a case may still be handled by a berthing maneuver. \\
As the rotation around the body fixed $\BFScoordinatey$--axis is supposed to be stable, we choose the values for the tensor of inertia $\inertiaxxtarget = \inertiazztarget$ and $\inertiaxxservicer = \inertiazzservicer$. Among others, these constants are shown in Table \ref{tb:constants}. 
\renewcommand{\arraystretch}{1.3}
\begin{table}[!htb]
	\begin{center}
	\caption{Constants of optimal control problem ($i \in \left\{T,S\right\}$)}\label{tb:constants}
	\begin{tabular}{ccl}
		variable & value & description \\\hline
		$a$ & 7071000 & orbit radius [m]\\
		$GM$ & $398 \cdot 10^{12}$ & gravitational constant [$N (m/kg)^2$]\\ 
		$n$ & $\sqrt{\frac{GM}{a^3}}$ & mean motion [1/s] \\
		$m$ & 100 & satellite mass [kg] \\
		$\controlposition_{max}$ & 0.1 & maximum thrust [N]\\
		$\controlmomentum_{max}$ & 1 & maximum torque [Nm] \\
		$J_{11}^{i}$ & 1000 & Angular mass around x [$kg/m^2$]\\
		$J_{11}^{i}$ & 2000 & Angular mass around y [$kg/m^2$]\\
		$J_{11}^{i}$ & 1000 & Angular mass around z [$kg/m^2$]\\
		$r^{i}$ & 1 & safety-area around target/servicer [m]\\
		$\boldsymbol{d}^{i}$ & $\left[0 \;\; 1.01 \;\; 0\right]^\top$ & docking point target/servicer [m] \\
	\end{tabular}
	\end{center}
\end{table}
Note that the values for $\inertiavector^{\target}$ and $\inertiavector^{\servicer}$ do not correspond to real satellites, but as only differences of them are of interest the assumption of having a rotational symmetric satellite, leads to a realistic behaviour. Due to simplicity we choose the safety areas around the satellites to be a sphere. Hence, the feasible set simplifies to
\begin{align*}
 \feasibleset = \left\{\state(t) \mid  \left\|\LVLHcoordinate(t)\right\| \geq 2 \right\} .
\end{align*}
Furthermore, we set the weights within the cost functional to $\weighttime=0$ and $\weightthrust = \weightmomentum = 1$ and add an additional maximal maneuver time of $7min$, i.e. $420s$. Since a maximal available time to perform a specific maneuver will be scheduled and one wants to minimize the consumed propellant, this is a reasonable assumption. We used a discretization of $N = 210$ intervals for the control, resulting in $1261$ optimization variables for the control and the final time.

\begin{figure}[!htb]
	\begin{center}
		\caption{Relative position and velocity of the docking points ($\LVLHcoordinatex$ -- dashed, $\LVLHcoordinatey$ -- dash-dotted, $\LVLHcoordinatez$ -- solid)}
		\label{fig:dockingrelative}
		\includegraphics[width=0.48\textwidth]{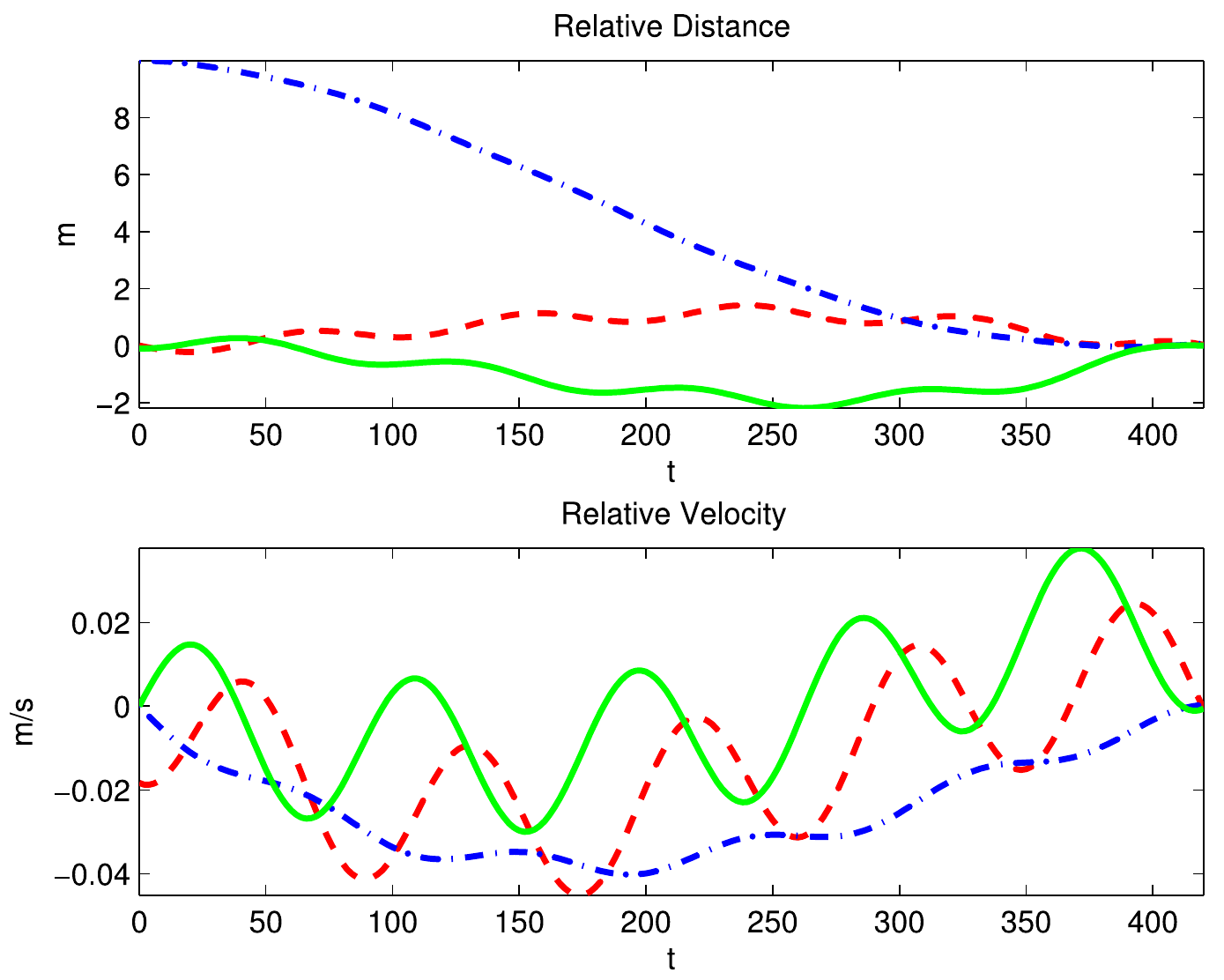}
	\end{center}
\end{figure}

\begin{figure}[!htb]
	\begin{center}
		\caption{Applied position and attitude control in body fixed coordinates ($\BFScoordinatex$ -- dashed, $\BFScoordinatey$ -- dash-dotted, $\BFScoordinatez$ -- solid)}
		\label{fig:control}
		\includegraphics[width=0.48\textwidth]{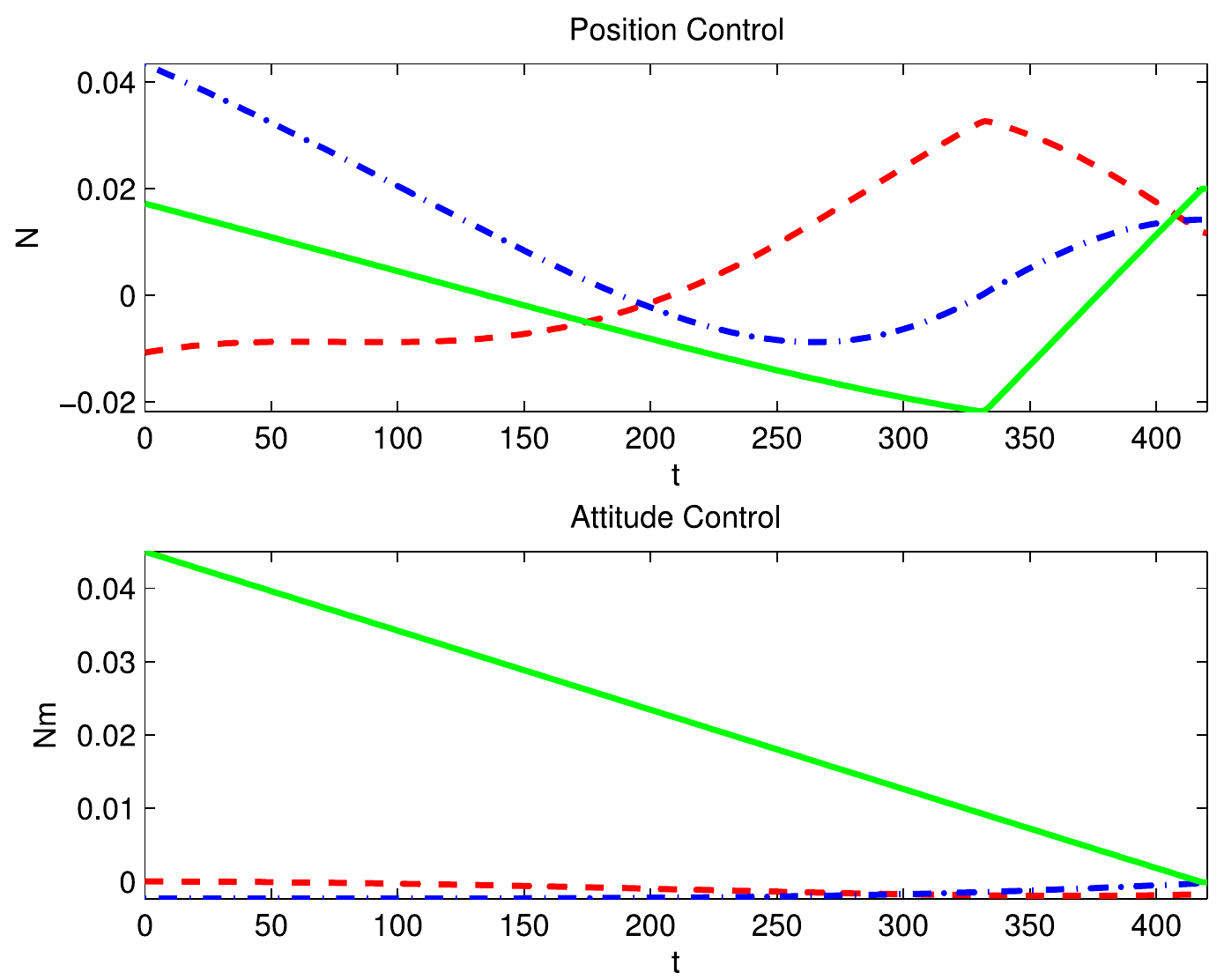}
	\end{center}
\end{figure}

\begin{figure*}[!htb]
	\begin{center}
		\caption{Trajectory of the servicer with docking vectors at $\finaltime$ and safety areas during transition}
		\label{fig:3dposition}
		\includegraphics[width=0.9\textwidth, height=6cm]{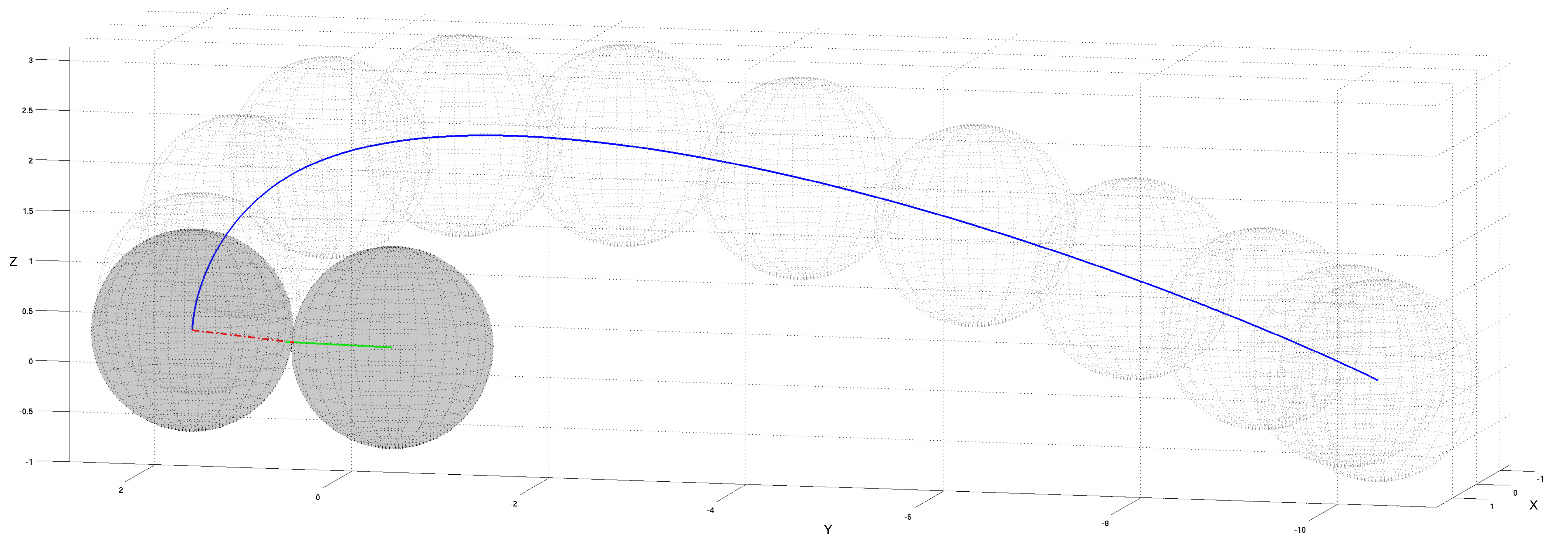}
	\end{center}
\end{figure*}

For the presented initial conditions, the computed results satisfy the Karush--Kuhn--Tucker conditions up to violations of magnitude $10^{-6}$ and the constraint violations are smaller than $10^{-8}$, cf. Figures \ref{fig:dockingrelative} -- \ref{fig:3dposition}. \\
Since the servicer has to surround the target to reach its destination, the state constraint influences the calculated trajectory in this example. Additional constraints concerning the concurrence of the orientation of the satellites cannot be applied in this case, as the servicer and the target have their dockingpoints located at their respective front. Due to the definition of the docking points lying only slightly outside the safety areas the possible angle between the docking vectors is limited. \\
Starting with relative distance and velocity trajectories of the docking maneuver, Figure \ref{fig:dockingrelative} shows that the relative distance in $\LVLHcoordinatey$--direction is continuously decreasing, whereas the $\LVLHcoordinatex$ and $\LVLHcoordinatez$ component ascend in the first half of the maneuver and then tend torwards the final value $0$. This behaviour is due to the dynamics of the Hill--Clohessy--Wilshire--Equations and the circumnavigation of the servicer around the target. 
Although the out of orbit motion $\LVLHcoordinatez$ is decoupled from the $\LVLHcoordinatex$--$\LVLHcoordinatey$ dynamic, satisfaction of the docking conditions upon termination is still required. Hence, not only positions but also velocities of the docking points have to coincide at $\finaltime$. Since the docking point of the target is tumbling, this results in the displayed oscillating behaviour in the relative position and velocity. The gyroscopic equations together with the above chosen values for the inertia tensor result in a periodic exchange of angular velocity between $\orientationxtarget$ and $\orientationztarget$.

The calculated controls in body fixed coordinates are shown in Figure \ref{fig:control}. Neither the position controls nor the attitude control are  on the boundary of the allowed control space. They prescribe a continuous and almost everywhere smooth control. This control structure seems to be realizable on a real satellite if the nozzle thrust can be adjusted continuously. Due to the chosen weights the consumed energy for position and attitude control is minimized and the terminal time is at the upper bound of the allowed maneuver time. The attitude control around the $\BFScoordinatez$ axis provides the main reorientation of the servicer, i.e. the turnaround of the satellite to ensure that both docking points face each other upon termination of the maneuver.

Last, Figure \ref{fig:3dposition} shows the trajectory prescribed by the servicer together with the terminal state of both satellites and some positions occupied by the servicer during the maneuver. One can see that the trajectory prescribes a smooth curve around the target and that the state constraint is satisfied during the approach. The two shaded spheres represent the bodies of the satellites in final position together with their current docking vectors. The solid and the dash dotted line represent the docking vectors $\dockingpoint^{\target}$ and $\dockingpoint^{\servicer}$ at the end of the maneuver. One can see that their endpoints coincide at time instant $\finaltime = 420$s.

\begin{figure}[!htb]
	\begin{center}
		\caption{Comparison of $\BFScoordinatex$--position--control for different control grids (dashed -- 50, dash--dotted -- 210, solid -- 420)}
		\label{fig:controlcomparison}
		\includegraphics[width=0.48\textwidth,height=4cm]{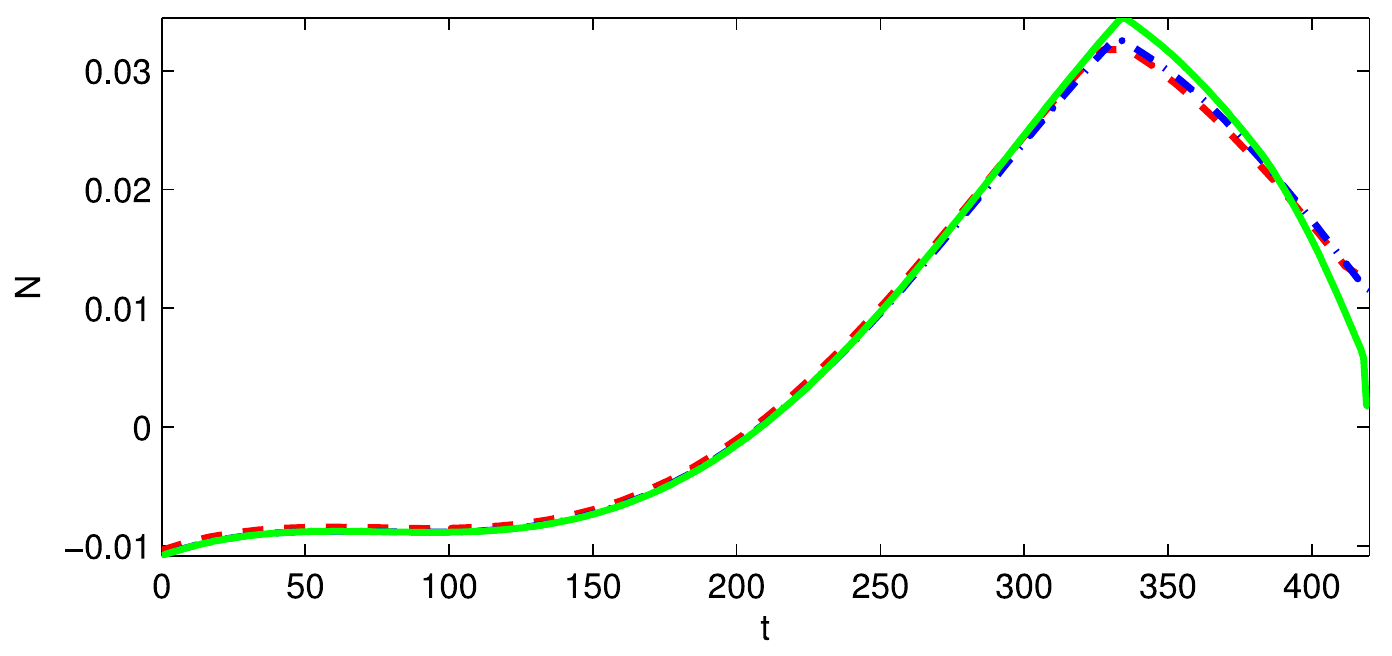}
	\end{center}
\end{figure}

The  computing times for the above shown example regarding refinements of the discretization can be seen in Table \ref{tb:computationtimes}. As expected, the required computing time increases in a nonlinear way with respect to the number of grid points. Additionally, one can see that even for coarse grids this method is not applicable to perform real time calculations, i.e. the control has to be calculated a priori. Yet, a comparison of the resulting controls shows that they differ only marginally as shown exemplarily for the position controls of the $\BFScoordinatex$--coordinate in Figure \ref{fig:controlcomparison}.
\begin{table}[!htb]
	\begin{center}
	\caption{Computation times for different control discretizations}\label{tb:computationtimes}
	\begin{tabular}{c|r}
	 N & computation time \\
         \hline
         50 & $1min$ $45sec$ \\
	210 & $2h$ $52min$ $9sec$ \\
        420 & $7h$ $44min$ $26sec$ \\
	\end{tabular}
   \end{center}
\end{table}

\section{Outlook}
\label{Section:outlook}
Future work will consider improving details of the geometry of the satellites. Until now collision test of complex geometries based on the union of convec polyhedreal sets can only be performed as post optimization verification. In future research, we plan to integrate this aspect into the optimization process. While incorporating these sets will enlarge the optimization problem drastically, the problem itself can be preprepared using Farka's Lemma to sort out inactive constraints. Additionally, improved initial guesses for the controls and thus the state trajectory may be imposed to enhance reliability of the algorithm and to reduce the computing time.

\bibliographystyle{ifacconf}
\bibliography{satellite-tumbling}

\end{document}